\documentclass [a4paper, 11pt]{article}

\usepackage{amsmath}
\usepackage{amsfonts}
\usepackage{amssymb}
\usepackage{amsthm}
\usepackage{makeidx}
\usepackage{graphicx}
\usepackage{pb-diagram}
\usepackage{epstopdf}
\usepackage{fancyhdr}
\usepackage{enumitem}
\usepackage[all]{xy}

\newtheorem{teo}{Theorem}[section]%definizioni che danno luogo ai nomi di teoremi, proposizioni, ecc.
\newtheorem{prop}[teo]{Proposition} %la quadra [teo] indica che la numerazione  subordinata a quella dei teoremi

\newtheorem{cor}[teo]{Corollary}
\theoremstyle{definition}
\newtheorem{defin}[teo]{Definition}

\newtheorem{nt}[teo]{Notation}

\newtheorem{nota}[teo]{Note}

\newtheorem{prob}[teo]{Problem}

\newtheorem{rem}[teo]{Remark}

   {\begin{verse}%
     \small }%
   {\end{verse}}

\providecommand{\U}{\mathcal{U}} \providecommand{\V}{\mathcal{V}}

\begin{document}

\setlength{\parskip}{1ex plus 0.5ex minus 0.2ex}
\begin{center}
$\mathbf{PRODUCT\,\,\,BETWEEN\,\,\,ULTRAFILTERS\,\,\,AND\,\,\,APPLICATIONS}$
$\mathbf{TO\,\,\,THE\,\,\,CONNES'\,\,\,EMBEDDING\,\,\,PROBLEM}$
\end{center}
\begin{center}
\end{center}
\begin{center}
V. CAPRARO - L. P\u AUNESCU\footnote{Second author is supported by the Marie Curie Research Training Network MRTN-CT-2006-031962 EU-NCG.}
\end{center}
$\mathbf{Abstract.}$ In this paper we want to apply the notion of
product between ultrafilters to answer several questions which arise
around the Connes' embedding problem. For instance, we will give a
simplification and generalization of a theorem by R\u adulescu; we
will prove that ultraproduct of hyperlinear groups is still
hyperlinear and consequently the von Neumann algebra of the free
group with uncountable many generators is embeddable into
$R^{\omega}$. This follows also from a general construction that
allows, starting from an hyperlinear group, to find a family of
hyperlinear groups. We will introduce the notion of hyperlinear pair
and we will use it to give some other characterizations of
hyperlinearity. We shall prove also that the cross product of a
hyperlinear group via a profinite action is embeddable into
$R^{\omega}$.

\section{Preliminaries}
We start by introducing the notion of product between ultrafilters.
It is already known in Model Theory (see, for example,
\cite{DiNa-Fo}), but it seems nobody applied it to Operator
Algebras.
\begin{defin}
Let $\mathcal U,\mathcal V$ be two ultrafilters respectively on $I$
and $J$. The tensor product $\mathcal U\otimes\mathcal V$ is the
ultrafilter on $I\times J$ defined by setting
$$
X\in\mathcal U\otimes\mathcal V\Leftrightarrow \{i\in I:\{j\in
J:(i,j)\in X\}\in\mathcal V\}\in\mathcal U
$$
\end{defin}
Observe that this is indeed a maximal filter, i.e. an ultrafilter.
\begin{rem}
This definition is equivalent to the following one:
$$
X\in\mathcal U\otimes\mathcal V\Leftrightarrow \exists A\in\mathcal
U\,\,\,s.\,\,t.\,\,\, \forall i\in A,
\pi_J(X\cap\pi_I^{-1}(i))\in\mathcal V
$$
where $\pi_I,\pi_J$ are the projections of $I\times J$ on the first
and
second component.\\
We prefer this second definition since it is easier to apply to
prove the following
\end{rem}
\begin{teo}\label{product}
Let $\{x_i^j\}_{(i,j)\in I\times J}\subseteq\mathbb R$ bounded. Then
$$
lim_{i\rightarrow\mathcal U}lim_{j\rightarrow\mathcal
V}x_i^j=lim_{(i,j)\rightarrow\U\otimes\V}x_i^j
$$
\begin{proof}
Let $x=lim_{i\rightarrow\U}lim_{j\rightarrow\V}x_i^j$. Fixed
$\varepsilon>0$, we notice from the definitions that
$$
A=\{i\in
I:|lim_{j\rightarrow\V}x_i^j-x|<\frac{\varepsilon}{2}\}\in\U
$$
and
$$
A_i=\{j\in
J|x_i^j-lim_{j\rightarrow\V}x_i^j|<\frac{\varepsilon}{2}\}\in\V
$$
Combining this two (by triangle inequality) we get
$$
X=\{(i,j)\in I\times J:i\in A,j\in A_i\}\subseteq\{(i,j)\in I\times
J:|x_i^j-x|<\varepsilon\}
$$
Since $X\in\U\otimes\V$ and $\varepsilon$ was arbitrary, it follows
the thesis.
\end{proof}
\end{teo}

\begin{nt}
By $\omega,\omega'$ we shall denote free ultrafilters on $\mathbb
N$. $R$ stands for the hyperfinite type $I I_1$ factor.  We shall
use the classical notation $R^\omega$ for the ultrapower of $R$ with
regard to $\omega$ and denote by $\tau$ its trace. By $L(G)$ we denote the von Neumann group
algebra of $G$.
\end{nt}

\section{Main result and immediate consequences}
The main result is actually an easy consequence of Th.\ref{product},
but it gives a tool to pass by the limit on representations. We
shall give some applications of this procedure.
\begin{prop}\label{ultraproduct}
Let $\omega,\omega'$ two ultrafilters on $\mathbb N$. Then
$$
(R^{\omega})^{\omega'}\cong R^{\omega\otimes\omega'}
$$
\begin{proof}
Those von Neumann algebras have the same algebraic structure. So we
only have to prove that they have the same trace. It is just a
consequence of Th.\ref{product}.
\end{proof}
\end{prop}

We want to apply this result to hyperlinear groups. In order to fully benefit
from it we will introduce the notion of hyperlinear pair.

\begin{defin}
By \emph{a central pair} we mean $(G,\varphi)$, where $G$ is a group
and $\varphi:G\to\mathbb C$ is a positive defined function, central
(i.e. constant on conjugacy classes) and $\varphi(e)=1$. Let
$Cen(G)$ be the set of those functions on $G$.
\end{defin}
\begin{rem}
An important element of $Cen(G)$ is the function $\delta_e$, defined
by setting $\delta_e(g)=0, \forall g\neq e$.
\end{rem}
\begin{rem}
If $(G,\varphi)$ is a central pair then we have a canonical
bi-invariant and bounded metric induced on $G$ by:
\[d(g,h)^2=2-\varphi(g^{-1}h)-\varphi(h^{-1}g)\, \,\,\,\,\,\,\,\,\,\, \forall g,h\in G.\]
\end{rem}

We recall that one can define the notion of ultraproduct of groups with bi-invariant metric (see \cite{Pe}).
We can use this definition for our particular case of central pairs.

\begin{defin}
Let ${(G_n,\varphi_n)}_{n\in\mathbb N}$ a sequence of central pairs and $\omega$ an
ultrafilter. By the \emph{ultraproduct} of the family we mean the central pair:
\[(G,\varphi)=(\Pi_{n}G_n/N,lim_\omega\varphi_n),\]
where $\Pi_{n}G_n$ is just the cartesian product and
$N=\{(g_n)\in\Pi_{n}G_n:lim_\omega\varphi_n(g_n)\to 1\}$.
\end{defin}
We shall denote by $\Pi_\omega(G_n,\varphi_n)$ the ultraproduct of central pairs.
\begin{nota}
It is easy to recognize that our definition of $N$ coincides with
the classical one: $N=\{(g_n)_n\in\Pi
G_n:lim_{\omega}d_n(g_n,e_n)\rightarrow0\}$.
\end{nota}
\begin{defin}
A central pair $(G,\varphi)$ is called \emph{hyperlinear} if there
exists an homomorphism $\theta_{\varphi}:G\rightarrow U(R^{\omega})$
such that
$$
\tau(\theta_{\varphi}(g))=\varphi(g)\,\,\,\,\,\,\,\,\,\,\,\,\,\forall
g\in G.
$$
Let $Hyp(G)=\{\phi\in
Cen(G):(G,\phi)\,\,\,is\,\,\,a\,\,\,hyperlinear\,\,\,pair\}$
\end{defin}
\begin{rem}
We recall the original definition by R\u adulescu: a countable
i.c.c. group $G$ is called hyperlinear if there exists a
monomorphism $G\rightarrow U(R^{\omega})$. It happens if and only if
$\delta_e\in Hyp(G)$ (see \cite{Ra}, Prop.2.5). Countability and
i.c.c. properties are not necessary, but they come from the reason
of this definition: to study when the group algebra is embeddable
into $R^{\omega}$. This problem, well-known as Connes' embedding
problem for groups, regard only separable type $II_1$ group factor.
\end{rem}
\begin{rem}
If $(G,\varphi)$ is a hyperlinear pair, then $(G,\varphi)$ is also a
central pair and the induced distance is just the distance in norm 2
in $R^\omega$.
\end{rem}

We can now use Prop. \ref{ultraproduct} in order to get the following

\begin{prop}\label{hyperlinear pairs}
Ultraproduct of hyperlinear pairs is a hyperlinear pair.
\begin{proof}
Take a sequence $(G_n,\varphi_n)$ of hyperlinear pairs and just embed each pair in an
$R^\omega$. The ultraproduct of the family with respect to $\omega'$ will sit inside
$(R^{\omega})^{\omega'}\cong R^{\omega\otimes\omega'}$.

In case we cannot find a "good" $\omega$ for all hyperlinear pairs, we just need to adapt
our notion of product between two ultrafilters to a notion of ultraproduct of ultrafilters. We
shall not do this, as it is just a technical trick and assuming continuum hypothesis this
$R^\omega$ are isomorphic between themselves anyway.
\end{proof}
\end{prop}
In order to give some information on the structure of $Hyp(G)$, we
recall that a \emph{monoid} is a set with a binary associative
operation admitting a neutral element. If $(X,\cdot)$ is a monoid,
an element $x\in X$ is called \emph{annihilator} if $x\cdot y=y\cdot
x=x, \forall y\in X$. The set of annihilators of $X$ is denoted by
$0(X)$. Clearly $Cen(G)$ is a monoid with respect the pointwise
product and $\delta_e\in0(Cen(G))$.

\begin{prop}\label{hyp(g)}
$Hyp(G)$ is a submonoid of $Cen(G)$. It is closed under ultralimits
and convex combinations. Moreover, for $G$ countable,
$0(Hyp(G))=\{\delta_e\}$ if and only if $G$ is hyperlinear in the
classical sense of R\u adulescu.
\begin{proof}
The constant function $1$ forms with $G$ a hyperlinear pair via the trivial representation.
$Hyp(G)$ is closed under pointwise multiplication because $R^\omega\otimes R^\omega\subset R^\omega$,
$\tau(x\otimes y)=\tau(x)\tau(y)$ and so $\theta_{\varphi\cdot\psi}=\theta_\varphi\otimes\theta_\psi$ will do the work.

For the second part note that $(G,lim_\omega\varphi_n)\subset\Pi_\omega(G,\varphi_n)$ and
use our last proposition. For convex combination define an homomorphism of $G$ in $R^\omega\oplus R^\omega$ with the same convex combination of traces.

The last part is an easy consequence of the following
Prop.\ref{annihilator}.
\end{proof}
\end{prop}

\begin{cor}
An i.c.c. group $G$ embeds in $U(R^\omega)$ if and only if $L(G)$ embeds
into $R^{\omega}$.
\begin{proof}
If $L(G)\subseteq R^{\omega}$ then clearly $G\subset U(R^\omega)$.
Conversely, let $\theta:G\rightarrow U(R^{\omega})$ an embedding.
Let $\tau$ be the normalized trace on $R^{\omega}$. Then
$|\varphi(g)|=|\tau(\theta(g))|<1$ for any $g\neq e$ (since $G$ is
i.c.c.) and $\varphi\in Hyp(G)$. Because of Prop. \ref{hyp(g)} we
have that $\varphi^n\in Hyp(G)$ and $lim_{n\to\omega}\varphi^n\in
Hyp(G)$.

Now $|\varphi(g)|<1$ for $g\neq e$ so $lim_n\varphi(g)^n=0$. This
means that $lim_n\varphi^n=\delta_e$, so $\delta_e\in Hyp(G)$. This is
equivalent to $L(G)$ embeds in $R^\omega$.
\end{proof}
\end{cor}

This is a simplification of the initial proof given by R\u adulescu in \cite{Ra} and also note that our proof doesn't need the contability of $G$.

\begin{prop}\label{annihilator}
A countable group $G$ is hyperlinear if and only if for any $g\in
G\setminus\{e\}$ there is a hyperlinear pair $(G,\varphi_g)$ such
that $|\varphi_g(g)|< 1$.
\begin{proof}
The \emph{only if} part is trivial. Conversely, we need to show that
$\delta_e\in Hyp(G)$. Take $G=\bigcup_n F_n$, with $F_n$ increasing
sequence of finite subsets of $G$. Define $\varphi_{F_n}=\Pi_{g\in
F_n}\varphi_g$. According to Prop.\ref{hyp(g)} $\varphi_{F_n}\in
Hyp(G)$ and by the same proposition so is
$\varphi=lim_{n\to\omega}\varphi_{F_n}$.

Now because of the hypothesis $|\varphi_g(g)|< 1$ and because of
$F_n$ is an increasing sequence we deduce $|\varphi(g)|< 1$. As in
the above corollary we now have $\delta_e=lim_n\varphi^n$, so
$\delta_e\in Hyp(G)$.
\end{proof}
\end{prop}

We end this section by presenting a motivation for our definition of
$Hyp(G)$. Let $\mathbb F_\infty$ be the free group with countable
many generators.

\begin{prop}
If $Cen(\mathbb F_\infty)=Hyp(\mathbb F_\infty)$ then every
countable group is hyperlinear.
\begin{proof}
Let $G$ be a countable group. Let $H$ be a normal subgroup of
$\mathbb F_\infty$ such that $G\cong\mathbb F_\infty/H$. Let
$\varphi_H:\mathbb F_\infty\to\mathbb C$ be the characteristic
function of $H$. We shall prove that $\varphi_H\in Cen(\mathbb
F_\infty)$. It is easy to see that $\delta_e\in Hyp(G)$ if and only
if $\varphi_H\in Hyp(\mathbb F_\infty)$. This will finish the proof.

Now $H$ is normal in $\mathbb F_\infty$. So for any $g,h\in \mathbb
F_\infty$ $h\in H$ if and only if $ghg^{-1}\in H$. This prove that
$\varphi_H$ is central. To prove that it is also positive defined
take $g_1,\ldots,g_n\in\mathbb F_\infty$. Consider the matrix
$\{\varphi_H(g_i^{-1}g_j)\}_{i,j}$ and notice that is the matrix of
an equivalence relation on a set with $n$ elements (because $H$ is a
subgroup). By permuting elements $(g_i)_i$ we can assume that is a
block matrix. This means that
$\sum_{i,j=1}^n\overline{\lambda_i}\lambda_j\varphi(g_i^{-1}g_j)$ is
nonnegative. So $\varphi_H$ is positive defined.
\end{proof}
\end{prop}

\begin{nota}
Our sets $Cen(G)$ and $Hyp(G)$ can be generalized to a type $II_1$
factor instead of just group algebras. Let $M$ be such a factor and
consider $B=\{x_n\}_{n\in\mathbb N}\subset M$ a basis in
$L^2(M,tr)$. Suppose that $x_0=id$. We shall consider now
$\varphi:B\to\mathbb C$ such that $\varphi(x_0)=1$ and the linear
extension of $\varphi$ to $M$ is positive and tracial (may not be
faithful). The problem is that such a linear extension may not be
well defined. We formalize this as follows: $\varphi\in Cen(M)$ iff
whenever $\varphi(x^*x)$ is well defined then so is $\varphi(xx^*)$
and $\varphi(x^*x)=\varphi(xx^*)\geq 0$.

For $\varphi\in Cen(M)$ we can define $M_\varphi$ by the $GNS$-construction. We define $\varphi\in Hyp(M)$ iff this $M_\varphi$ is embedable in $R^\omega$. As we saw, for $M=L(G)$ and $\varphi_H$ for $H$ a normal subgroup of $G$ then $L(G)_{\varphi_H}=L(G/H)$.

As another example we may take the crossed product $M=L^\infty(X)\rtimes G$ of a non-free measure preserving action. Take $\{f_i:i\in\mathbb N\}$ a basis for $L^\infty(X)$ and $B=\{f_iu_g:i\in\mathbb N,\ g\in G\}$. Define $\varphi(f_iu_g)=\int_{X_g}f_i$ where $X_g=\{x\in X:gx=x\}$. Then $M_\varphi=M(E_G)$, the Feldmann-Moore construction for the equivalence relation induced by $G$ on $X$.
\end{nota}

\section{Other applications}

\subsection{Construction of uncountable hyperlinear groups}
Now we want to present a construction that, starting from an
hyperlinear group $G$, allows to construct a family of countable and
uncountable hyperlinear groups. An easy application of this
construction is that the von Neumann algebra of the free group with
uncountable many generators $\mathbb F_{\aleph_c}$ is embeddable
into $R^{\omega}$. The Hilbert-Schmidt distance between two distinct
universal unitaries of $\mathbb F_{\aleph_c}$ will be equal to
$\sqrt{2}$, giving another proof of the non-separability of $R^{\omega}$.

\begin{defin}
Let $G$ be a countable group with generators $g_1,g_2,...$. Let
$\Im$ be a family of infinite subsets of $\mathbb N$ such that
$F_1,F_2\in\Im$ implies $F_1\cap F_2$ is finite. Now let
$F=\{f_1,f_2,...\}\in\Im$, define the sequence $(g_n^F)_n=g_{f_n}$.
Let $g^F$ be the sequence $g_n^F$ modulo $\omega$. We can multiply
$g^{F_1}$, $g^{F_2}$ component-wise, by using the relations on $G$.
The group generated by the elements $g^F$ is denoted by $G(\omega,\Im)$.\\\\
Notice that $G(\omega,\Im)$ does not depend only on $\omega$ and
$\Im$, but also on the set of generators chosen.
\end{defin}

\begin{rem}
The generators $g^F$ of $G(\omega,\Im)$ are different elements in
$G(\omega,\Im)$. This is because $g_n^{F_1}=g_n^{F_2}$ holds only
for a finite number of indexes, by the definition of $\Im$. Since a
free ultrafilter does not contain finite sets, $g^{F_1}$ and
$g^{F_2}$ must be different.
\end{rem}

\begin{rem}
$G(\omega,\Im)$ can be countable (if the family $\Im$ is countable),
but also uncountable. Indeed one can use the Zorn's lemma to prove
the existence of an uncountable family $\Im$ which verifies the
property $F_1,F_2\in\Im$ implies $F_1\cap F_2$ is finite. An elegant
example privately suggested by Ozawa is the following: take
$t\in[\frac{1}{10},1)$, for example $t=0,132483...$, define
$$
I_t=\{1,13,132,1324,13248,132483,...\}
$$
i.e. $I_t$ is the set of the approximation of $t$. Then
$\{I_t\}_{t\in[\frac{1}{10},1)}$ is an uncountable family of subsets
of $\mathbb N$ such that $I_t\cap I_s$ is finite for all $t\neq s$.
\end{rem}

\begin{prop}
If $G$ is hyperlinear, then also $G(\omega,\Im)$ is hyperlinear.
\begin{proof}
We want to prove that $G(\omega,\Im)\subset\Pi_\omega(G,\delta)$ and
the last is a hyperlinear pair because of Prop.\ref{hyperlinear
pairs}. Moreover we shall prove that if in an ultraproduct of
central pairs just $\delta_e$ appears, then the central positive
defined function of the ultraproduct will also be $\delta_e$. This two
affirmations will show that $\delta_e\in Hyp(G(\omega,\Im))$, i.e.
$G(\omega,\Im)$ is hyperlinear.

Recall that
$\Pi_\omega(G_n,\varphi_n)=(\Pi_{n}G_n/N,lim_\omega\varphi_n)$,
where $\Pi_{n}G_n$ is just the cartesian product and
$N=\{(g_n)\in\Pi_{n}G_n:lim_\omega\varphi_n(g_n)\to 1\}$. So let
$G_n$ a copy of $G$ and $\varphi_n=\delta_e$ for each $n$. Then
$lim_\omega\varphi_n\in\{0,1\}$. If this limit is $1$ for some
element, then that element is in $N$ i.e. it is the identity in the
ultraproduct. So indeed $lim_\omega\delta_e=\delta_e$ proving our second
affirmation.

Now from the construction of $G(\omega,\Im)$ we see that
$G(\omega,\Im)\subset \Pi_nG_n$. If an element $g=(g_n)_n$ of
$G(\omega,\Im)$ is in $N$ then $lim_\omega\delta_e(g_n)=1$ meaning
that $g_n=e$ in $G$ for any $n$ in a set in $\omega$. From the
definition of $G(\omega,\Im)$ this means that $g=e$. We proved that
$G(\omega,\Im)\subset\Pi_\omega(G,\delta)$.
\end{proof}
\end{prop}

It is well known that $\mathbb F_\infty$, free group with countable
many generators is hyperlinear. We shall denote with $\mathbb
F_{\aleph_c}$ the free group with $\aleph_c$ many generators (set of
continuum power).

\begin{cor}
$\mathbb F_{\aleph_c}$ is hyperlinear. In particular $R^{\omega}$ is
not separable.
\begin{proof}
If $Card(\Im)=\aleph_c$ then $\mathbb F_\infty(\omega,\Im)=\mathbb F_{\aleph_c}$,
and we can apply the previous proposition.

Representing $L(\mathbb F_{\aleph_c})$ on $R^\omega$, the
Hilbert-Schmidt distance between two elements of $\mathbb
F_{\aleph_c}$ will be $\sqrt{2}$. Separability in the weak or in the
strong topology is the same and the last one coincide with the
Hilbert-Schmidt topology on the bounded sets (see \cite{Jo}).
\end{proof}
\end{cor}

\begin{nota}
Non-separability of $R^{\omega}$ is already well-known. The first proof is probably due to Feldman (see \cite{Fel}); S. Popa proved in \cite{Po} that every MASA in $R^{\omega}$ is not separable. Anyway, we want to underline the importance of non-separability of $R^\omega$ around the Connes' embedding conjecture: every separable type $II_1$ factor can be embedded into $R^{\omega}$ (see \cite{Co}). This conjecture imply the existence of a universal type $II_1$ factor. If a factor embeds in $R^\omega$ then it embeds in any $R^{\omega'}$. We are grateful to Pestov for communicating this fact to us. Ozawa proved in \cite{Oz} that such a universal factor cannot be separable, also proved by Nicoara, Popa and Sasyk in \cite{Ni-Po-Sa}).  So, if $R^{\omega}$ was been separable, Connes embedding conjecture would be false.
\end{nota}

\begin{prob}
What kind of groups have the shape $G(\omega,\Im)$? Is it true that
if $\{R_a\}_{a\in A}$ is the set of distinct relations on $G$ and
$B\subseteq A$, then there exist $\omega$ and $\Im$ such that the
set of relations of $G(\omega,\Im)$ is $\{R_a\}_{a\in B}$?
\end{prob}

\subsection{Cross product via profinite actions}

We want to apply Prop.\ref{ultraproduct} also to some other type $II_1$ factors
than group algebras. For this we ask ourselves when the crossed product
$L^{\infty}(X)\rtimes_{\alpha}G$ for a free action $\alpha$ embeds in $R^\omega$.
Of course when this happens $G$ has to be hyperlinear. We shall prove the converse
in the easy case in which $\alpha$ is profinite.

\begin{defin}
Let $\alpha$ be an action of a group $G$ on a von Neumann algebra
$P$. Then $\alpha$ is called \emph{profinite} if there is an
increasing sequence of finite dimensional $G$-invariant subalgebras
$A_1\subset A_2\subset\ldots$ such that $P=(\bigcup_n A_n)''$.
\end{defin}

\begin{prop}\label{profinite}
Let $G$ be a hyperlinear group and $\alpha$ be a profinite action of
$G$ on $X$. Then $L^{\infty}(X)\rtimes_{\alpha}G$ is embeddable into
$R^{\omega}$.
\begin{proof}
The crossed product is generated on $L^2(X)\otimes l^2G$ by the
operators $\alpha(g)\otimes\lambda(g)$ for $g\in G$ and $m_f\otimes
1$ for $f\in L^\infty(X)$ (here $\lambda$ is the regular
representation of $G$ on $l^2G$ and $m_f$ is the multiplication
operator).

Let $L^\infty(X)=(\bigcup_nA_n)''$ with $A_n$ G-invariant and finite
dimensional. We can then form $A_n\rtimes_{\alpha} G$ and
$L^\infty(X)\rtimes_{\alpha} G=(\bigcup_n A_n\rtimes_{\alpha} G)''$.
Looking at the above definition of crossed product we can deduce
that $A_n\rtimes_{\alpha} G \subset M_{k_n}\otimes L(G)$. Here
entered the fact that $A_n$ is finite dimensional. Now, because $G$
is hyperlinear $M_{k_n}\otimes L(G)\subset R\otimes R^\omega\subset
R^\omega$. We can than embed $\bigcup_n A_n\rtimes_{\alpha} G$ in
$(R^\omega)^{\omega'}$ so that $L^\infty(X)\rtimes_{\alpha} G\subset
R^{\omega\otimes\omega'}$.
\end{proof}
\end{prop}

\section{Acknowledgements}

It is our pleasure to thanks professor Florin R\u adulescu for many
useful discussions and remarks on the topics presented in this
paper.

V. CAPRARO, \emph{UNIVERSIT\`{A} DI ROMA TOR VERGATA} e-mail:
capraro@mat.uniroma2.it

L. P\u AUNESCU, \emph{UNIVERSIT\`{A} DI ROMA TOR VERGATA and INSTITUTE of MATHEMATICS "S. Stoilow" of the ROMANIAN ACADEMY} email:
paunescu@mat.uniroma2.it


\begin{thebibliography}{9}
\bibitem[Co]{Co} A. Connes, \emph{Classification of injective factors}, Ann. of Math. 104
(1976), 73-115.
\bibitem[DiNa-Fo]{DiNa-Fo} M. Di Nasso - M. Forti, \emph{Hausdorff ultrafilters},
Proc. Amer. Math. Soc.  134  (2006), 1809-1818.
\bibitem[Fa-Ha-Sh]{Fa-Ha-Sh} I. Farah - B.Hart - D. Sherman, \emph{Model theory of operator algebras I:
Stability}, arXiv:math/0908.2790
\bibitem[Fe]{Fe} J. Feldman, C. Moore, \emph{Ergodic equivalence Relations, Cohomology, and Von Neumann Algebras II},
Trans. Amer. Math. Soc. Vol. 234, No. 2 (1977) pp.325-359.
\bibitem[Fel]{Fel} J. Feldman, \emph{Nonseparability of certain finite
factors}, Proc. Amer. Math. Soc. 7 (1956), 23--26.
\bibitem[Ge-Ha]{Ge-Ha} L. Ge - D. Hadwin, \emph{Ultraproducts of
$C^*$-algebras}, Oper. Theory Adv. Appl. 127 (2001), 305-326.
\bibitem[Io]{Io} A. Ioana, \emph{Cocyle Superrigidity for Profinite Actions of property (T) Groups}, arxiv:0805.2998 (2008)
\bibitem[Jo]{Jo} V.R. Jones, \emph{von Neumann algebras}, notes from
a course.
\bibitem[Ni-Po-Sa]{Ni-Po-Sa}, R. Nicoara - S. Popa - R. Sasyk,
\emph{on type $II_1$ factors arising from 2-cocycles of w-rigid
groups}, J. Funct. Anal. 242 (2007), no.1, 230-246.
\bibitem[Oz]{Oz} N. Ozawa, \emph{There is no separable $II_1$
universal factor}, Proc. Amer. Math. Soc. 132 (2) (2004),
arXiv:math/0210411v2.
\bibitem[Pe]{Pe} V. Pestov, \emph{Hyperlinear and Sofic Groups: A Brief Guide}, arXiv:math/0804.3968v8(2008).
\bibitem[Po]{Po} S. Popa, \emph{On a problem of R. V. Kadison on maximal abelian *-subalgebras in
factors}, Invent. Math. 65 (1981/82), 269–281.
\bibitem[R\u a]{Ra} F. Radulescu, \emph{The von Neumann algebras of the non-residually finite Baumslag group $<a,b|ab^3a^{-1}=b^2>$ embeds into
$R^{\omega}$}, arXiv:math/0004172v3 (2000).
\end{thebibliography}
\end{document}